\documentstyle[amsfonts,11pt]{article}
\newcommand{\ad}{{\mathrm a}{\mathrm d\,}}

 \makeatother

\newtheorem{thm}{Theorem}[section]

\newtheorem{prop}[thm]{Proposition}

\begin{document}
\title{The classification of $\omega$-Lie algebras}

\author{Zhiqi Chen  \\ School of Mathematics and Statistics, \\ Guangdong University of Technology, Guangzhou 510520, P.R. China. \\ E-mail: chenzhiqi@nankai.edu.cn \\
Junna Ni \\ Department of Mathematics, South China Normal University, \\ Guangzhou 510631, P.R.China. Corresponding author. Email: nijunna@126.com  \\ Jianhua Yu \\ Department of Mathematics, South China Normal University, \\ Guangzhou 510631, P.R.China. Email: yujianhuscnu@126.com}

\maketitle

\begin{abstract}
Zusmanovich gave a fundamental result on the structure of $\omega$-Lie algebras. But up to now, the classification of $\omega$-Lie algebras is still open. In this paper, we give a complete classification of $\omega$-Lie algebras over $\mathbb C$.


\end{abstract}

\section{Introduction}
A vector space $L$ over $\mathbb F$ is called an $\omega$-Lie algebra if there is a bilinear map $[\cdot,\cdot]:L\times L\rightarrow L$ and a
skew-symmetric bilinear form $\omega: L \times L \rightarrow \mathbb F$ such that
\begin{enumerate}
  \item $[x,y]=-[y,x]$,
  \item $[[x,y],z]+[[y,z],x]+[[z,x],y]=\omega(x,y)z+\omega(y,z)x+\omega(z,x)y$,
\end{enumerate}
hold for any $x,y,z\in L$. Clearly, an $\omega$-Lie algebra $L$ with $\omega=0$ is a Lie algebra, which is called a trivial $\omega$-Lie algebra. Otherwise, $L$ is called a nontrivial $\omega$-Lie algebra.

The notation of an $\omega$-Lie algebra is given by Nurowski in \cite{Nur}, and there are many studies in this field such as \cite{chen1, 4dim,  chen3, zhang2, Zu}. In particular, \cite{Zu} plays a fundamental role on the study of the structure of $\omega$-Lie algebras. Also it is pointed out in \cite{Zu} that it is difficult to give the classification of nontrivial $\omega$-Lie algebras. In this paper, the classification of $\omega$-Lie algebras is given as follows.

\begin{thm}\label{classification}
A nontrivial $\omega$-Lie algebra $L$ over $\mathbb C$ is one of the following cases:
\begin{enumerate}
  \item $\dim L=3$, $L$ is $L_1$, or $L_2$, or $A_{\alpha}$, or $B$, or $C_{\alpha}$ in Proposition~\ref{3dim}.
  \item $\dim L=4$, $L$ is $L_{1,i}(i=1,2,\cdots,8)$, or $L_{2,i}(i=1,2,3,4)$, or $E_{1,\alpha}(\alpha\not=0,1)$, or $F_{1,\alpha}(\alpha\not=0,1)$, or $G_{1,\alpha}$, or $H_{1,\alpha}$, $\widetilde{A_{\alpha}}$, or $\widetilde{B}$, or $\widetilde{C_{\alpha}}(\alpha\not=0,-1)$ in Proposition~\ref{4dim}.
  \item  $\dim L\geq 5$, $L=H\oplus {\mathbb C}x\oplus {\mathbb C}v$, and the nonzero brackets and $\omega$ are:  $$[x,v]=h_1+x,\ [x,h]=f(h),\ [v,h]=g(h),\ \omega(x,v)=1.$$
       Here $h_1\in H$ and $f,g: H \rightarrow H$. Let $A$ and $B$ be the matrices decided by $f$ and $g$. Then $A,B\in {\mathfrak gl}(H)$, $A=I+A'$ and $[A',B]=A'$, where $A'$ is strictly upper triangular and $B$ is upper triangular. Denote this $\omega$-Lie algebra by $NP_1$-type.
  \item $\dim L\geq 5$, $L={\mathbb C}h_0\oplus H_1\oplus {\mathbb C}x\oplus {\mathbb C}v$, and the nonzero brackets and $\omega$ are:
       $$[x,h_0]=-ah_0, \ [v,h]=\frac{1}{a}h,\ [v,h_0]=h_2+\frac{1}{a}h_0+x, \ [x,v]=h_1+av,\ \omega(x,v)=1,$$ for any $h\in H_1$. Here $h_1\in {\mathbb C}h_0\oplus H_1$, $h_2\in H_1$ and $0\not=a\in {\mathbb C}$. Denote this $\omega$-Lie algebra by $P_1$-type.
  \item $\dim L\geq 5$, and $L=H\oplus {\mathbb C}x\oplus {\mathbb C}y\oplus{\mathbb C}a$, and the nonzero brackets and $\omega$ are: \begin{eqnarray*}
       && [a,h]=h, \ [x,h]=f(h),\ [y,h]=g(h), \ [y,a]=h_2,\\
       && [x,y]=h_3+d_1x-b_2d_1y+a,\ [x,a]=h_1-x+b_2y,\ \omega(x,y)=1,
      \end{eqnarray*}
     where $h_1,h_2,h_3\in H$, $f,g: H \rightarrow H$. Denote this $\omega$-Lie algebra by $NP_2$-type.
  \item $\dim L\geq 5$, and $L=H\oplus {\mathbb C}x\oplus {\mathbb C}y\oplus{\mathbb C}a$, and the nonzero brackets and $\omega$ are: \begin{eqnarray*}
       &&[a,h]=h,\ [x,h]=f(h),\ [y,a]=h_2-y,\\
       &&[x,a]=h_1+b_2y,\ [x,y]=h_3+d_2y+a,\ \omega(x,y)=1,
        \end{eqnarray*}
      where $h_1,h_2,h_3\in H$, $f: H \rightarrow H$ and $f(h_2)=d_2h_2.$ Denote this $\omega$-Lie algebra by $NP_3$-type.
  \item $\dim L\geq 5$, and $L=H\oplus {\mathbb C}x\oplus {\mathbb C}y\oplus{\mathbb C}a$, and the nonzero brackets and $\omega$ are:
        $$[a,h]=h, \ [x,y]=h_3+a,\ [x,a]=h_1+b_1x+b_2y,\ [y,a]=h_2+c_1y,\ \omega(x,y)=1,$$
       for any $h\in H$. Here $h_1,h_2,h_3\in H$, $b_1\not=0$, $c_1\not=0$ and $b_1+c_1+1=0$. Denote this $\omega$-Lie algebra by $P_2$-type.
  \item $\dim L\geq 5$, and $L=H\oplus {\mathbb C}x\oplus {\mathbb C}y\oplus{\mathbb C}a$, and the nonzero brackets and $\omega$ are:   \begin{eqnarray*}
          &&[a,h]=h,\ [x,h]=\lambda(h)y+f(h), \ [y,a]=-y,\\
          && [x,a]=c_2y+h_4,\ [x,y]=d_1y+a+h_3, \ \omega(x,y)=1,
       \end{eqnarray*}
    for any $h\in H$. Here $h_3,h_4\in H$, $\lambda: H\rightarrow {\mathbb C}$ is not zero, $f: H\rightarrow H$. Denote this $\omega$-Lie algebra by $NP_4$-type.
\end{enumerate}
\end{thm}

Throughout this paper, all vector spaces and algebras are finite dimensional over $\mathbb C$ unless stated otherwise.

\section{Preliminary}

The classification of $\omega$-Lie algebras of dimensions $1$ and $2$ is obvious. For the 3-dimensional case, we have the following proposition.

\begin{prop}[\cite{chen1}]\label{3dim}
Let $L$ be a nontrivial 3-dimensional $\omega$-Lie algebras. Then it is isomorphic to one of the following cases: there exists a basis $\{x,y,z\}$ of $L$ such that the nonzero brackets and $\omega$ are given as follows.
\begin{enumerate}
  \item $L_1: [y,z]=z,[x,y]=y,\omega(x,y)=1.$
  \item $L_2: [x,z]=y, [y,z]=z, \omega(x,z)=1.$
  \item $A_{\alpha}: [x, y] = x, [x,z] = x + y, [y,z] = z + \alpha x, \omega(y,z) = -1$, where $\alpha \in \mathbb C$.
  \item $B: [x, y] = y, [x,z] = y + z, [y,z] = x, \omega(y,z) =2.$
  \item $C_{\alpha}: [x, y] = y, [x,z] = \alpha z, [y,z] = x, \omega(y,z) = 1 + \alpha$, where $0, -1\not=\alpha\in \mathbb C$.
 \end{enumerate}
\end{prop}

Otherwise, for $\omega$-Lie algebras of dimension 4,
\begin{prop}[\cite{4dim}]\label{4dim}
Let $L$ be a 4-dimensional nontrivial $\omega$-Lie algebras. Then it is isomorphic to one of the following algebras: there exists a basis $\{x,y,z,e\}$ of $L$ such that the nonzero brackets and $\omega$ are
\begin{enumerate}
 \item $L_{1,1}: [x,y]=y,[y,z]=z,[e,y]=-y,\omega(x,y)=1.$
  \item $L_{1,2}: [x,y]=y, [y,z]=z, [e,x]=z,[e,y]=-e,\omega(x,y)=1.$
  \item $L_{1,3}: [x,y]=y,  [y,z]=z,  [e,x]=y,  [e,y]=-e,\omega(x,y)=1,  \omega(e,x)=1.$
  \item $L_{1,4}:  [x,y]=y,  [y,z]=z,  [e,x]=y+z,  [e,y]=-e, \omega(x,y)=1, \omega(e,x)=1.$
  \item $L_{1,5}: [x,y]=y, [y,z]=z,  [e,x]=e,  [e,y]=-e, \omega(x,y)=1$
  \item $L_{1,6}:  [x,y]=y,  [y,z]=z, [e,x]=e+y,  [e,y]=-e,\omega(x,y)=1, \omega(e,x)=1.$
  \item $L_{1,7}: [x,y]=y,  [y,z]=z,  [e,x]=e,  [e,y]=z-e,  \omega(x,y)=1.$
  \item $L_{1,8}: [x,y]=y,  [y,z]=z,  [e,x]=e+y,  [e,y]=z-e, \omega(x,y)=1,  \omega(e,x)=1.$
  \item $L_{2,1}: [x,z]=y,  [y,z]=z,  [e,y]=-e,  \omega(x,z)=1.$
  \item $L_{2,2}: [x,z]=y,  [y,z]=z,  [e,y]=-e, [e,x]=z,  \omega(x,z)=1.$
  \item $L_{2,3}: [x,z]=y,  [y,z]=z,  [e,y]=-e,  [e,x]=e, \omega(x,z)=1.$
  \item $L_{2,4}: [x,z]=y, [y,z]=z,  [e,y]=-e,  [e,x]=e+z, \omega(x,z)=1.$
  \item $E_{1,\alpha}(\alpha\neq 0,1): [x,y]=y,  [y,z]=z,  [e,x]=\alpha e,  [e,y]=-e, \omega(x,y)=1.$
  \item $F_{1,\alpha}, (\alpha\neq 0,1): [x,y]=y,  [y,z]=z,  [e,y]=-e,  [e,x]=\alpha e+y,\omega(x,z)=1, \omega(e,x)=1.$
  \item $G_{1,\alpha}: [e, x]=e+\alpha y, [e, y]=-e+x, [y,z] = z, [x, y] = y, \omega(e, x)=\alpha, \omega(x, y) = 1.$
  \item $H_{1,\alpha}: [e, x] = e + \alpha y,[e, y]=-e + x + z, [y,z] = z, [x, y] = y,  \omega(e, x) = \alpha, \omega(x, y) = 1.$
  \item $\widetilde{A_{\alpha}}: [x,y]=x, [x,z]=x+y, [y,z]=z+{\alpha}x,  [e,z]=e, \omega(y,z)=-1$.
  \item $\widetilde{B}: [x,y]=y, [x,z]=y+z, [y,z]=x, [e,x]=-2e, [e,y]=-e, \omega(y,z)=2.$
  \item $\widetilde{C_{\alpha}}\ ({\alpha}\neq 0,-1): [x,y]=y, [x,z]={\alpha}z, [y,z]=x, [e,x]=-(1+{\alpha})e, \omega(y,z)=1+{\alpha}.$
\end{enumerate}
\end{prop}

The fundamental result on nontrivial $\omega$-Lie algebras is described in \cite{Zu}.
\begin{prop}[\cite{Zu}]\label{2.1}
 Let $L$ be a finite-dimensional nontrivial $\omega$-Lie algebra. Then
one of the following holds:
\begin{enumerate}
  \item $\dim L = 3.$
  \item $L$ has a Lie subalgebra of codimension 1.
  \item $\ker \omega$ is an almost abelian Lie algebra of codimension 2 in $L$ with the abelian part acting nilpotently on $L$.
\end{enumerate}
\end{prop}

Moreover, then Lie subalgebra of codimension 1 in Case 2 of Proposition~\ref{2.1} must be one of the following cases (\cite{Zu}):
\begin{enumerate}
   \item $L' = H \oplus {\mathbb{C}}x$, $H$ is abelian, $\ad x : H \rightarrow H$ is any linear map, and $\ker\omega=H$.
   \item $L'$ is the direct sum of an abelian Lie algebra $H$ and the two-dimensional non-abelian Lie algebra $\langle x, y | [x, y] = y\rangle$, and $\ker\omega=H\oplus {\mathbb{C}}x$.
   \item $L' = H \oplus {\mathbb{C}}x\oplus{\mathbb{C}}y$, $H$ is abelian, $\ad x : H \rightarrow H$ is the identity map, $\ad y : H \rightarrow H$
is any linear map, $[x, y] \in H$, and $\ker\omega=H\oplus {\mathbb{C}}x$.
   \item $L' = H \oplus {\mathbb{C}}x\oplus{\mathbb{C}}y$, $H$ is abelian, $\ad x : H \rightarrow H$ is the identity map, $\ad y : H \rightarrow H$
is the zero map, $[x, y] = h_0+ky$ for some $h_0 \in H, 0\not=k\in {\mathbb C}$, and $\ker\omega=H\oplus {\mathbb{C}}x$.
\end{enumerate}

Further study gives
\begin{prop}[\cite{Zu}]
   A finite dimensional semisimple $\omega$-Lie algebra is either a Lie algebra, or has dimension $\leq 4$.
\end{prop}

But it is also pointed out in~\cite{Zu} that it is far from giving the complete classification of $\omega$-Lie algebras.

\section{The classification of nontrivial $\omega$-Lie algebras}
In order to classify $\omega$-Lie algebras, it is enough to classify $\omega$-Lie algebras of dimension $\geq 5$ since the classification of $\omega$-Lie algebras for $\dim \leq 4$ is clear. Let $L$ be an $\omega$-Lie algebra with $\dim L\geq 5$.

\subsection{Case 2 in Proposition~\ref{2.1}.}
First for {\bf cases 3 and 4}, $\ker\omega=H\oplus {\mathbb{C}}x$ is an almost abelian Lie algebra of codimension 2, which is a solvable Lie algebra. Clearly $L$ is a representation of $\ker\omega$ under the adjoint action. For any $h\in H$, $h=[x,h]$, therefore $\ad h$ is nilpotent on $L$ by Lie Theorem. Thus, these two cases include in case 3 in Proposition~\ref{2.1}.

For {\bf case 1},  $L = H \oplus {\mathbb{C}}x \oplus {\mathbb{C}}v$ with $\ker\omega=H$ and $\omega(x,v)=1.$ For any $h\in H$, let $$[x,h]=f(h),\quad [v,h]=g(h)+\lambda(h)x+\mu(h)v.$$ Here $f,g: H\rightarrow H$ and $\lambda, \mu: H\rightarrow {\mathbb C}.$

$1)$ If there exists $h_0\in H$ such that $\mu(h_0)=1,$ then $H=H_1\oplus{\mathbb C}h_0$, where $\mu(h)=0$ for any $h\in H_1.$
That is, for any $h\in H_1$, $$[v,h_0]=h'+ax+v,\quad [v,h]=g(h)+\lambda(h)x, \quad [x,v]=h_1+b_1x+b_2v.$$
Let $x'=x+b_2h_0$. Then we have
\begin{eqnarray*}
  && [x',h]=f(h),\forall h\in H,\quad [v,h_0]=h''+ax'+v, \\
  && [v,h]=g_1(h)+\lambda(h)x',\forall h\in H_1,\quad [x',v]=h_2+bx'.
\end{eqnarray*}
Here $g_1: H\rightarrow H$. By the $\omega$-Jacobi identity, for $h\in H_1$, we have
\begin{eqnarray*}
0 & = & [[v,h_0],h]+[[h_0,h],v]+[[h,v],h_0] =[h''+ax'+v,h]+[-g_1(h)-\lambda(h)x',h_0]\\
  & = & af(h)+g_1(h)+\lambda(h)x'-\lambda(h)f(h_0).
\end{eqnarray*}
It follows that $\lambda(h)=0$ for any $h\in H_1$, i.e. $[v,h]=g_1(h)$ for any $h\in H_1.$ Let $f(h_0)=h_3+th_0$ for some $h_3\in H_1.$ Then
\begin{eqnarray*}
h_0 & = & [[x',v],h_0]+[[v,h_0],x']+[[h_0,x'],v]\\
&  = & bf(h_0)+[h'',x']+[v,x']+[-h_3-th_0,v]\\
& = & b(h_3+th_0)-f(h'')-h_2-bx+g_1(h_3)+t(h''+ax'+v).
\end{eqnarray*}
It follows that $t=0$, and then $b=0$ which means that $L$ is not perfect. That is, we have that $L=H_1\oplus{\mathbb C}h_0\oplus {\mathbb{C}}x \oplus {\mathbb{C}}v$, and for any $h\in H_1$, \begin{eqnarray*}
&&[v,h_0]=h_1+ax+v,\quad[x,v]=h_2,\quad [x,h_0]=h_3\in H_1,\\
&&[x,h]=f(h),\quad [v,h]=g(h),\quad \omega(x,v)=1.
\end{eqnarray*}Furthermore, for any $h\in H_1,$
\begin{eqnarray*}
0&=&[[v,h_0],h]+[[h_0,h],v]+[[h,v],h_0]\\
&=&[h_1+ax+v,h]+[-g(h),h_0]\\
&=&af(h)+g(h).
\end{eqnarray*}
Replacing $v$ by $ax+v,$ we can assume $[v,h_0]=h_1+v.$
Then for any $h\in H_1$, $$[v,h]=0$$ by $af(h)+g(h)=0.$ Clearly for any $h\in H_1$, $$[x,h]=f(h),\quad [x,v]=h_2,\quad [x,h_0]=h_3,\quad\omega(x,v)=1.$$
For any $h\in H_1$, let $f(h)=f_1(h)+t_hh_0,$ where $f_1(h)\in H_1$ and $t_h\in\mathbb C$, then
\begin{eqnarray*}
h&=&[[x,v],h]+[[v,h],x]+[[h,x],v]=[[h,x],v]\\
&=&[-f(h),v]=-[f_1(h)+t_hh_0,v]\\
&=&t_h(h_1+v).
\end{eqnarray*}
It follows that $t_h=0$ and then $h=0$ which is impossible.

%师弟第二页最后，有问题。后续又修正如下：
$2)$ For any $h\in H$, $\mu(h)=0$. That is, $[v,h]=g(h)+\lambda(h)x$ for any $h\in H.$ If $\lambda(h)=0$ for any $h\in H$, then we have that 
$L=H\oplus {\mathbb C}x\oplus {\mathbb C}v$, and for any $h\in H$
 \begin{eqnarray*}
    &&[x,h]=f(h),\quad[v,h]=g(h),\quad [x,v]=h_1+ax+bv,\quad \omega(x,v)=1.
 \end{eqnarray*}
For this case, $L$ is not perfect. Then for any $h\in H$, we have that
\begin{eqnarray*}
h&=&[[x,v],h]+[[v,h],x]+[[h,x],v]\\
&=&[h_1+ax+bv,h]+[g(h),x]+[-f(h),v]\\
&=&af(h)+bg(h)-f(g(h))+g(f(h)).
\end{eqnarray*}
It is impossible when $a=b=0$. We can suppose $a\neq 0$ by the symmetry of $x$ and $V$.
Replacing $x$ and $v$ by $ax+bv$ and $\frac{1}{a}v$ respectively, we can assume that
  \begin{eqnarray*}
&&[x,v]=h_1+x,\quad[x,h]=f(h),\quad [v,h]=g(h),\quad \omega(x,v)=1.
 \end{eqnarray*}
Here we also denote $f'$ and $g'$ by $f$ and $g$. Then the $\omega$-Lie Jacobi identity for $x,v,h$ is
$$f(h)-f(g(h))+g(f(h))=h.$$
Let $A$ and $B$ be the matrixes of $\ad x$ and $\ad v$ on $H$ respectively. Then the above equation is $A-AB+BA=id,$ namely
$$[A-id,B]=A-id.$$ That is, $A-id, B$ span a solvable  Lie subalgebra of $gl(H)$. Thus $A-id, B$ can be upper triangular relative to a suitable basis of $H$. Then $A-id=[A-id,B]$ is a strictly upper triangular matrix denoted by $A'$. Then we have that
 \begin{eqnarray*}
&&A=id+A',\quad[A',B]=A',
 \end{eqnarray*}
where $A'$ is strictly upper triangular and $B$ is upper triangular. {\bf This is $NP_1$-type in Theorem~\ref{classification}.}

Otherwise there exists $h_0\in H$ such that $[v,h_0]=h_1+x$ for some $h_1\in H.$ Furthermore, there exists $H_1$ such that $H=H_1\oplus {\mathbb C}h_0$, and
$$[v,h]=g(h), \forall h\in H_1. $$
Here $g: H_1\rightarrow H$. Assume that $$[x,h]=f(h),\quad [x,v]=h'_2+b_1x+b_2v,$$ where $f: H\rightarrow H.$
Take $x'=x+b_1h_0,$ then we get that
\begin{eqnarray*}
&&[x',v]=h_2+b_2v, \quad [v,h_0]=h_3+x',\quad [x',h]=f(h),\quad \omega(x',v)=1.
\end{eqnarray*}
Assume that $f(h_0)=h_4+th_0$ for some $h_4\in H_1.$ Then
\begin{eqnarray*}
 -h_0&=&[[x',h_0],v]+[[h_0,v],x']+[[v,x'],h_0]\\
 &=&[f(h_0),v]+[-h_3-x',x']+[-h_2-b_2v,h_0]\\
 &=&-g(h_4)-t(h_3+x')+f(h_3)-b_2(h_3+x').
\end{eqnarray*}
Then $t=-b_2,$ i.e. $f(h_0)=h_4-b_2h_0.$ Then for any $h\in H_1,$
\begin{eqnarray*}
 0&=&[[h_0,v],h]+[[v,h],h_0]+[[h,h_0],v]\\
 &=&[-h_3-x',h]+[g(h),h_0]\\
 &=&-f(h).
\end{eqnarray*}
That is, $[x',h_0]=h_4-b_2h_0$ and $[x',h]=0$ for any $h\in H_1.$ 
If $b_2=0$, then for any $h\in H_1$, we have
$$h=[[x',v],h]+[[v,h],x']+[[h,x'],v]=[[v,h],x']\in {\mathbb C}h_4$$
which is impossible for $\dim H_1\geq 2$, i.e. $\dim L\geq 5.$
Thus we must have $b_2\neq 0.$
Take $h'_0=h_0-b_2^{-1}h_4$. Then $H=H_1\oplus {\mathbb C}h_0'$, $[x',h_0']=-b_2h_0'$, and
$[v,h_0']=h_3'+x'$ for some $h_3'\in H$. Assume that $h_3'=h_3''+th_0'$ for some $h_3''\in H_1$. Then
\begin{eqnarray*}
 -h_0'&=&[[x',h_0'],v]+[[h_0',v],x']+[[v,x'],h_0']\\
 &=&[-b_2h_0',v]+[-h_3''-th_0'-x',x']+[-h_2-b_2v,h_0']\\
 &=&-tb_2h_0',
\end{eqnarray*}
then $t=b_2^{-1}$. That is, $[v,h_0']=h_3''+b_2^{-1}h_0'+x'$ for some $h_3''\in H_1.$ For any $h\in H_1$, let $g(h)=\phi(h)+\lambda(h)h_0',$ where $\phi: H_1\rightarrow H_1, \lambda: H_1\rightarrow {\mathbb C}.$
\begin{eqnarray*}
 -h&=&[[h,v],x']+[[v,x'],h]+[[x',h],v]\\
 &=&[-\phi(h)-\lambda(h)h_0',x']+[-h_2-b_2v,h]\\
&=&-\lambda(h)b_2h_0'-b_2[v,h]\\
&=&-\lambda(h)b_2h_0'-b_2(\phi(h)+\lambda(h)h_0')\\
&=&-b_2\phi(h)-2\lambda(h)b_2h_0'.
\end{eqnarray*}
Then $\lambda(h)=0$ and $\phi(h)=b_2^{-1}h$ if $\dim H\geq 2$. For this case, $L$ is perfect. In summary, $H=H_1\oplus {\mathbb C}h_0'$, and for some $b_2\not=0$ and for any $h\in H_1$
\begin{eqnarray*}
 && [H,H]=0, \quad [x',h_0']=-b_2h_0', \quad [x',h]=0, \\
 && [v,h]=b_2^{-1}h,\quad [v,h_0']=h_3''+b_2^{-1}h_0'+x', \quad [x',v]=h_2+b_2v.
\end{eqnarray*}
{\bf It is $P_1$-type in Theorem~\ref{classification}.}

For {\bf case 2}, $L = H \oplus {\mathbb{C}}x\oplus {\mathbb{C}}y\oplus {\mathbb C}v$, and we have $$[x,h]=0,\quad [y,h]=0,\quad [x,y]=y,\quad \omega(y,v)=1.$$

If there exists $h_0\in H$ such that $[h_0,v]=h_1+ax+by+v$, then we can assume
\begin{eqnarray*}
&&[h,v]=f(h)+\lambda(h)x+\mu(h)y, h\in H_1,
\end{eqnarray*}
where $H=H_1\oplus {\mathbb C}h_0$, $f: H_1\rightarrow H$, and $\lambda, \mu: H_1\rightarrow {\mathbb C}.$ Replace $x$ by $x-eh_0$ for some $e$ if necessary, we can assume $$[x,v]=h_2+cx+dy.$$
For any $h\in H,$
\begin{eqnarray*}
0 & = & [[v,x],h]+[[x,h],v]+[[h,v],x]\\
  & = & [-h_2-cx-dy,h]+[[h,v],x]\\
  & = & [[h,v],x].
\end{eqnarray*}
In particular,
$0=[[h_0,v],x]=[h_1+ax+by+v,x]=-by-h_2-cx-dy.$
It gives that
$ c=0$, $b=-d$, and $h_2=0$. That is, $$[x,v]=-by.$$
%For any $h\in H_1$,
%$0=[[h,v],x]=[f(h)+\lambda(h)x+\mu(h)y,x]=-\mu(h)y$, so $\mu(h)=0$. That is, for any $h\in H_1,$ $$[h,v]=f(h)+\lambda(h)x.$$
Assume that $[v,y]=h_3+tx+ky+mv.$
%For any $h\in H_1$, we have
%\begin{eqnarray*}
% -h & = & [[v,y],h]+[[y,h],v]+[[h,v],y] \\
%    & = & [[v,y],h]+[f(h)+\lambda(h)x,y] \\
%    & = & [[v,y],h]+\lambda(h)y
%\end{eqnarray*}
%thus $\lambda(h)=0, h\in H_1.$ That is $[h,v]=f(h), h\in H_1.$
Then we have
\begin{eqnarray*}
-x&=&[[v,y],x]+[[y,x],v]+[[x,v],y]\\
  &=&[h_3+tx+ky+mv,x]+[-y,v]+[-by,y]\\
  &=&-ky+mby+h_3+tx+ky+mv.
\end{eqnarray*}
It gives that $m=0$, $t=-1$, and $h_3=0$. That is, $$[v,y]=-x+ky.$$
It implies that
\begin{eqnarray*}
-h_0&=&[[v,y],h_0]+[[y,h_0],v]+[[h_0,v],y]\\
  &=&[-x+ky,h_0]+[h_1+ax+by+v,y]\\
  &=&ay-x+ky,
\end{eqnarray*}
which is a contradiction.

Then the projection of $[h,v]$ on ${\mathbb C}v $ is zero for any $h\in H.$ Namely, we can assume that for any $h\in H$,
$$[h,v]=f(h)+\lambda(h)x+\mu(h)y.$$
%[H,v]在v上的投影为零，假设[x,v]|_v=0
Let $[y,v]=h_1+ax+by+ev$. If $[x,v]=h_2+cx+dy$, then by $\omega$-Jacobi identity,
\begin{eqnarray*}
-x&=&[[x,v],y]+[[v,y],x]+[[y,x],v]\\
  &=&[h_2+cx+dy,y]+[-h_1-ax-by-ev,x]-[y,v]\\
  &=&cy+by+e(h_2+cx+dy)-(h_1+ax+by+ev),
\end{eqnarray*}
we get $e=0.$ But, for any $h\in H$, $$h=[[y,v],h]+[[v,h],y]+[[h,y],v]=[[v,h],y]=\lambda(h)y,$$
which is impossible. So we can assume that
%[x,v]|_v\neq 0,
$$[x,v]=h_2+cx+dy+v,$$
By $\omega$-Jacobi identity,
\begin{eqnarray*}
x&=&[[y,v],x]+[[v,x],y]+[[x,y],v]\\
  &=&[h_1+ax+by+ev,x]+[-h_2-cx-dy-v,y]+[y,v]\\
  &=&(b+c)[y,x]+e[v,x]+2[y,v]\\
  &=&-(b+c)y-e(h_2+cx+dy+v)+2(h_1+ax+by+ev),
\end{eqnarray*}
we also have $e=0.$  For any $h\in H,$
\begin{eqnarray*}
0&=&[[x,v],h]+[[v,h],x]+[[h,x],v]\\
  &=&[h_2+cx+dy+v,h]+[-f(h)-\lambda(h)x-\mu(h)y,x]\\
  &=&-f(h)-\lambda(h)x,
\end{eqnarray*}
which means $\lambda(h)=0$ and $f(h)=0.$ Namely, $$[h,v]=\mu(h)y.$$
Again for any $h\in H,$
$$h=[[y,v],h]+[[v,h],y]+[[h,y],v]=[h_1+ax+by,h]+[-\mu(h)y,y]=0,$$
which is impossible.

\subsection{Case 3 in Proposition~\ref{2.1}}
For this case, by \cite{Zu}, we have that $L=H\oplus {\mathbb C}x\oplus{\mathbb C}y\oplus {\mathbb C}a$, and
$$\ker \omega=H\oplus {{\mathbb C}}a, [H,H]=0, [a,h]=h,\forall h\in H.$$
If $\dim L \geq 5$, we have the following cases:
\begin{enumerate}
  \item[I.] $[x,h]=f(h)$, $[y,h]=g(h)$.
  \item[II.] $[x,h]=\lambda(h)y+f(h)$, $[y,h]=g(h)$, $g(\ker \lambda)=0$.
\end{enumerate}
Here $\lambda\not=0: H \rightarrow  {\mathbb C} $ and $f,g: H\rightarrow H$.

{\bf Case I.} For any $h\in H$, by the $\omega$-Jacobi identity, we have
$$ [[y,a],h]=[y,[a,h]]+[[y,h],a]=[y,h]+[g(h),a]=0.$$
Similarly, $[[x,a],h]=0.$ Clearly $<H,a>$ is a solvable Lie subalgebra of $L$. By Lie Theorem, we can set
$$[x,a]=h_1+b_1x+b_2y+b_3a, \quad [y,a]=h_2+c_1y+c_2a. $$
Together with $[[y,a],h]=[[x,a],h]=0$, we have
\begin{eqnarray*}
 & c_1g(h)+c_2h=0,  & \quad(a)  \\
 & b_1f(h)+b_2g(h)+b_3h=0. & \quad(b)
 \end{eqnarray*}

If $c_1=0,$ by (a), $c_2=0,$ then $[L,L]\neq L$, i.e. $L$ is not perfect. For this case, $[y,a]=h_2.$ Let $[x,y]=h_3+d_1x+d_2y+d_3a.$
By the $\omega$-Jacobi identity,
\begin{eqnarray*}
[[x,a],y]&=&[x,[a,y]]+[a,[y,x]]+\omega(y,x)a\\
&=&[x,-h_2]+[h_3+d_1x+d_2y+d_3a,a]-a\\
&=&-f(h_2)-h_3+d_1[x,a]+d_2[y,a]-a\\
&=&-f(h_2)-h_3+d_1(h_1+b_1x+b_2y+b_3a)+d_2h_2-a\\
&=&-f(h_2)-h_3+d_1h_1+d_2h_2+d_1b_1x+d_1b_2y+(d_1b_3-1)a.
\end{eqnarray*}
On the other side,
\begin{eqnarray*}
[[x,a],y]&=&[h_1+b_1x+b_2y+b_3a,y]=-g(h_1)+b_1[x,y]-b_3h_2\\
         &=&-g(h_1)+b_1(h_3+d_1x+d_2y+d_3a)-b_3h_2\\
         &=&-g(h_1)+b_1h_3 -b_3h_2+b_1d_1x+b_1d_2y+b_1d_3a.
\end{eqnarray*}
It follows that $$b_1d_2=b_2d_1,\quad b_1d_3=b_3d_1-1.\quad (b)$$
If $b_1=0,$ then by Eq.(b),$$b_2=0,\quad d_1\neq 0.$$ Furthermore, $b_3=\frac{1}{d_1}$ .
That is $[x,a]=h_1+\frac{1}{d_1}a.$ But
\begin{eqnarray*}
0=[[x,a],h]=[h_1+\frac{1}{d_1}a,h]=\frac{1}{d_1}h,
\end{eqnarray*}
which is impossible. So $b_1\neq 0$. Then $$d_2=\frac{b_2d_1}{b_1}, \quad d_3=\frac{b_3d_1-1}{b_1}.$$
That is  $$[x,y]=h_3+d_1x+\frac{b_2d_1}{b_1}y+\frac{b_3d_1-1}{b_1}a.$$
Take $x'=x+b_1^{-1}b_3a$, $f_1=f+b_1^{-1}b_3$ and $h_3^{'}=h_3-b_1^{-1}b_3h_2.$ It is easy to check that
\begin{eqnarray*}
&&[x',a]=h_1+b_1x'+b_2y,\quad [x',h]=f_1(h),\\
&&[x',y]=h_3^{'}+d_1x'+\frac{b_2d_1}{b_1}y-\frac{1}{b_1}a,\quad \omega(x',y)=1.
\end{eqnarray*}
Furthermore, by
\begin{eqnarray*}
&&0=[[x',a]h]=[h_1+b_1x'+b_2y,h]=b_1f_1(h)+b_2g(h),
\end{eqnarray*}
 we get $f_1(h)=-b_1^{-1}b_2g(h).$
Then for $[[x',y],h]$, we have
\begin{eqnarray*}
[[x',y],h]&=&[h_3^{'}+d_1x'+\frac{b_2d_1}{b_1}y-\frac{1}{b_1}a,h]\\
          &=&d_1f_1(h)+\frac{b_2d_1}{b_1}g(h)-\frac{1}{b_1}h\\
          &=&d_1\frac{-b_2}{b_1}g(h)+\frac{b_2d_1}{b_1}g(h)-\frac{1}{b_1}h\\
          &=&-\frac{1}{b_1}h.
\end{eqnarray*}
By the $\omega$-Jacobi identity,
\begin{eqnarray*}
[[x',y],h]&=&[x',[y,h]]+[y,[h,x']]+h\\
          &=&[x',g(h)]+[y,-f_1(h)]+h\\
          &=&f_1(g(h))-g(f_1(h))+h\\
          &=&-\frac{b_2}{b_1}g(g(h))-g(-\frac{b_2}{b_1}g(h))+h\\
          &=&h.
\end{eqnarray*}
It follows that $$b_1=-1.$$
In summary, there exist $x,y,a\in L$ such that, for any $h\in H$,
\begin{eqnarray*}
&&[x,h]=f(h),\quad [y,h]=g(h), \quad [H,H]=0, \quad [a,h]=h, \quad [x,y]=h_3+d_1x-b_2d_1y+a,\\
&&[x,a]=h_1-x+b_2y,\quad[y,a]=h_2, \quad \ker \omega=H\oplus {{\mathbb C}}a, \quad\omega(x,y)=1,
\end{eqnarray*} 
where $h_1,h_2,h_3\in H$. {\bf This is $NP_2$-type in Theorem~\ref{classification}.}

If $c_1\neq 0,$ by (a), $$g(h)=-c_1^{-1}c_2 h.$$
If $b_1=0,$ by (b), $b_2c_2=c_1b_3$. It follows that $[L,L]\neq L$, i.e. $L$ is not perfect. And we have
$$[x,a]=h_1+b_2y+b_3a,\quad [y,a]=h_2+c_1y+c_2a.$$
From (a) and (b), we have that
\begin{eqnarray*}
&&g(h)=-c_1^{-1}c_2h,\quad b_2c_2=b_3c_1.
\end{eqnarray*}
Let $y'=y+c_1^{-1}c_2a.$ Then it follows
\begin{eqnarray*}
&&[x,a]=h_1+b_2(y+c_1^{-1}c_2a)+b_3a-b_2c_1^{-1}c_2a=h_1+b_2y',\\
&&[y',h]=[y,h]+c_1^{-1}c_2[a,h]=g(h)+c_1^{-1}c_2h=0,\\
&&[y',a]=h_2+c_1y',\quad \omega (x,y')=1. 
\end{eqnarray*}
Let $[x,y']=h_3+d_1x+d_2y'+d_3a.$  Also we have
\begin{eqnarray*}
[[x,y'],a]&=&[h_3+d_1x+d_2y'+d_3a,a]\\
          &=&-h_3+d_1[x,a]+d_2[y',a]\\
          &=&-h_3+d_1(h_1+b_2y')+d_2(h_2+c_1y')\\
          &=&-h_3+d_1h_1+d_2h_2+(d_1b_2+d_2c_1)y'.
\end{eqnarray*}
On the other hand,
\begin{eqnarray*}
[[x,y'],a]&=&[x,[y',a]]+[[x,a],y']+\omega(x,y')a\\
          &=&[x,h_2+c_1y']+[h_1+b_2y',y']+a\\
          &=&f(h_2)+c_1(h_3+d_1x+d_2y'+d_3a)+a\\
          &=&f(h_2)+c_1h_3+c_1d_1x+c_1d_2y'+(c_1d_3+1)a.
\end{eqnarray*}
It follows that
\begin{eqnarray*}
&&c_1d_1=0,\quad c_1d_2=d_1b_2+c_1d_2,\quad c_1d_3+1=0,\\
&&f(h_2)+c_1h_3=-h_3+d_1h_1+d_2h_2.
\end{eqnarray*}
Then $d_1=0$ and $d_3=-\frac{1}{c_1}$.
Furthermore,
\begin{eqnarray*}
[[x,y'],h]&=&[h_3+d_1x+d_2y'+d_3a,h]=d_1f(h)+d_3h\\
         &=&[x,[y',h]]+[[x,h],y']+\omega(x,y')h=h
\end{eqnarray*}
Thus $d_3h=h$ by $d_1=0.$ For $h\neq 0,$ we have that $d_3=1,c_1=-1$ and $f(h_2)=d_2h_2.$ In summary, we get that
\begin{eqnarray*}
&&[y',a]=h_2-y',\quad [x,a]=h_1+b_2y',\quad [y',h]=0,\\
&&[x,h]=f(h),\quad [a,h]=h,\quad [x,y']=h_3+d_2y'+a,
\end{eqnarray*}
for $f(h_2)=d_2h_2.$ {\bf This is $NP_3$-type in Theorem~\ref{classification}.}

For the case $b_1\neq 0,$ by (b), we have  $$f(h)=-b_1^{-1}b_2g(h)-b_1^{-1}b_3h.$$
Take $y'=y+c_1^{-1}c_2a$ and $x'=x+(b_1^{-1}b_3-b_1^{-1}b_2c_1^{-1}c_2)a$. It easy to check that
\begin{eqnarray*}
 &&[x',h]=[y',h]=0,\\
 && [x',a]=h_1+b_1x'+b_2y',\\
 && [y',a]=h_2+c_1y'.
\end{eqnarray*}
Assume that $[x',y']=h_3+d_1x'+d_2y'+d_3a$. Then $[[x',y'],h]=d_3h.$ On the other side,
$$[[x',y'],h]=[x',[y',h]]+[y'[h,x']]+\omega(x',y')h=h.$$ It means $d_3=1.$ Furthermore by $\omega$-Jacobi identity,
\begin{eqnarray*}
[[x',y'],a]&=&[x',[y',a]]+[y',[a,x']]+\omega(x',y')a\\
&=&[x',h_2+c_1y']+[h_1+b_1x'+b_2y',y']+a\\
&=&c_1[x',y']+b_1[x',y']+a\\
&=&(c_1+b_1)h_3+(c_1+b_1)d_1x'+(c_1+b_1)d_2y'+(c_1+b_1+1)a.
\end{eqnarray*}
On the other side,
\begin{eqnarray*}
[[x',y'],a]&=&[h_3+d_1x'+d_2y'+a,a]\\
&=&-h_3+d_1(h_1+b_1x'+b_2y')+d_2(h_2+c_1y')\\
&=&-h_3+d_1h_1+d_2h_2+d_1b_1x'+(d_1b_2+d_2c_1)y'.
\end{eqnarray*}
It follows that
$$(c_1+b_1)d_1=b_1d_1, \quad b_1d_2=b_2d_1, \quad c_1+b_1+1=0.$$
By the first identity, we have $c_1d_1=0$, then $d_1=0$ since $c_1\not=0$. It follows that $d_2=0$ since $b_1\not=0$. Then $[x',y']=h_3+a$.

That is, for this case, $L$ is perfect, and there exist $x,y,a\in L$ such that, for any $h\in H$,
\begin{eqnarray*}
&&[x,h]=0,\quad [y,h]=0, \quad [H,H]=0, \quad [a,h]=h, \quad [x,y]=h_3+a,\\
&&[x,a]=h_1+b_1x+b_2y,\quad[y,a]=h_2+c_1y, \quad \ker \omega=H\oplus {\mathbb{C}}a, \quad\omega(x,y)=1.
\end{eqnarray*}
Here $h_1,h_2,h_3\in H$, $b_1\not=0$, $c_1\not=0$ and $b_1+c_1+1=0$. {\bf This is $P_2$-type in Theorem~\ref{classification}.}

{\bf Case II.} For this case, we claim that $L$ is not perfect. In fact, first by $\omega$-Jacobi identity, we have $[[y,a],h]=0$ for any $h\in H$.
Let $[y,a]=h_1+b_1x+b_2y+b_3a$. Then
\begin{eqnarray*}
[[y,a],h]&=&[h_1+b_1x+b_2y+b_3a,h]\\
         &=&b_1(\lambda(h)y+f(h))+b_2g(h)+b_3h\\
         &=&b_1\lambda(h)y+b_1f(h)+b_2g(h)+b_3h.
\end{eqnarray*}
Since $\dim L\geq 5$ and $\lambda\not=0$, we can choice $h\in H$ such that $\lambda (h)\neq 0$. Then by $[[y,a],h]=0$ and the above identity, we have $$b_1=0, \quad b_2g(h)+b_3h=0.$$
By $\dim L\geq 5,$ $\dim H\geq 2$, and then $\dim\ker \lambda \geq 1.$  There must exists $0\not=h\in H$ such that $g(h)=0.$ Then $0=b_2g(h)+b_3h=b_3h$ means $b_3=0$. That is, $$[y,a]=h_1+b_2y.$$

{\bf 1. $b_2\neq 0,-1$.} By $[[y,a],h]=[h_1+b_2y,h]=b_2[y,h]=b_2g(h)=0,$ we get that $g(h)=0$, i.e. $[y,h]=0.$
Let $y'=y+\frac{1}{b_2+1}h_1$ and $f_1(h)=f(h)-\frac{\lambda(h)}{b_2+1}h_1$. Then it is easy to check that
\begin{eqnarray*}
&& [x,h] =  \lambda(h)y'+f_1(h), \quad  [y',h]=  [y,h]=0, \\
&& [y',a]=   b_2y', \quad  \omega(x,y')=  \omega(x,y)=1.
\end{eqnarray*}
 Then by the $\omega$-Jacobi identity, we have
 \begin{eqnarray*}
 [[x,y'],h]&=&[x,[y',h]]+[y',[h,x]]+\omega(x,y')h\\
&=&-[y,\lambda(h)y+f_1(h)]+h\\
&=&h.
\end{eqnarray*}
Let $[x,y']=d_1y'+d_2a+d_3x+h_3$. Then
 \begin{eqnarray*}
 [[x,y'],h]&=&[d_1y'+d_2a+d_3x+h_3,h]=d_3[x,h]+d_2h\\
           &=& d_3\lambda(h)y'+d_3f_1(h)+d_2h.
\end{eqnarray*}
It follows that $d_3=0$ and $d_2=1$. That is,
  $$[x,y']=h_3+d_1y'+a.$$
Let $[x,a]=h_4+c_1x+c_2y'+c_3a.$ Similarly we have
\begin{eqnarray*}
[[x,a],h]&=&[x,[a,h]]+[a,[h,x]]=\lambda(h)(1+b_2)y' \\
        &=&[h_4+c_1x+c_2y'+c_3a,h]=c_1\lambda(h)y'+c_1f_1(h)+c_3h.
\end{eqnarray*}
It follows that
$$1+b_2=c_1.$$
Then for $[[x,a],y']$, we have
 \begin{eqnarray*}
 [[x,a],y']&=&[h_4+(1+b_2)x+c_2y'+c_3a,y']\\
 &=&((1+b_2)d_1-b_2c_3)y'+(1+b_2)a+(1+b_2)h_3 \\
 &=&[x,[a,y']]+[a,[y',x]]+\omega(y',x)a\\
  &=&-b_2[x,y']+[d_1y'+a+h_3,a]-a\\
 &=&-b_2(d_1y'+a+h_3)+d_1b_2y'-h_3-a\\
  &=&-(1+b_2)a-(1+b_2)h_3.
 \end{eqnarray*}
Then we have that $1+b_2=0$, which is a contradiction.

{\bf 2. $b_2=-1.$} Namely
 $$[a,h]=h, \quad [x,h]=\lambda(h)y+f(h), \quad [y,a]=h_1-y.$$
Then $[y,h]=0$ by $[[y,a],h]=0$. Let $[x,y]=d_1y+d_2a+d_3x+h_3.$ For $[[x,y],h]$, we have
\begin{eqnarray*}
[[x,y],h] &=& d_2h+d_3\lambda(h)y+d_3f(h) \\
          &=&[x,[y,h]]+[y,[h,x]]+\omega(x,y)h\\
          &=&h.
\end{eqnarray*}
It follows that $d_3=0$ and $d_2=1.$ That is,
$$[x,y]=d_1y+a+h_3.$$
Let $[x,a]=c_1x+c_2y+c_3a+h_4$. Then we have

\begin{eqnarray*}
[[x,a],h]&=&[x,[a,h]]+[a,[h,x]] =\lambda(h)y+f(h)-\lambda(h)[a,y]-[a,f(h)]\\
&=&\lambda(h)y+f(h)+\lambda(h)h_1-\lambda(h)y-f(h) =\lambda(h)h_1 \\
&=&c_1[x,h]+c_2[y,h]+c_3[a,h] =c_1\lambda(h)y+c_1f(h)+c_3h.
\end{eqnarray*}
Then $c_1\lambda(h)=0$ which implies $c_1=0.$ Then for $0\not=h\in \ker \lambda$, $c_3h=0$, i.e. $c_3=0$. That is,
$$[x,a]=c_2y+h_4,$$ which means $[L,L]\neq L$, i.e. $L$ is not perfect. Furthermore $h_1=0$. In fact, we have
\begin{eqnarray*}
0=[[x,a],h]&=&[x,[a,h]]+[a,[h,x]]\\
         &=&\lambda(h)y+f(h)+[\lambda(h)y+f(h),a]\\
         &=&\lambda(h)y+f(h)+\lambda(h)(h_1-y)-f(h)\\
         &=&\lambda(h)h_1.
\end{eqnarray*}
Then $h_1=0$ by taking $h\in H$ such that $\lambda(h)\neq 0.$ That is, there exist $x,y,a$ such that
\begin{eqnarray*}
&&[y,a]=-y,\quad [x,a]=c_2y+h_4,\quad [x,y]=d_1y+a+h_3,\\
&&[a,h]=h,\quad [x,h]=\lambda(h)y+f(h),\quad [y,h]=0,\quad \omega(x,y)=1.
\end{eqnarray*}
{\bf This is $NP_2$-type in Theorem~\ref{classification}.}

{\bf 3. $b_2=0.$} Then we have
$$[x,h]=\lambda(h) y+f(h),\quad [y,h]=g(h),\quad [a,h]=h, \quad [y,a]=h_1, \quad g(\ker \lambda)=0.$$
Let $[x,a]=c_1x+c_2y+c_3a+h_4$ and $[x,y]=d_1y+d_2x+d_3a'+h_3$. For $h_0\in H$ satisfying $\lambda(h_0)=1$, let $a'=a-c_2h_0$ and $y'=y-d_1h_0$. Then it is easy to check that
\begin{eqnarray*}
&& [x,a']= c_1x+c_3a'+h'_4\ (h'_4=c_3c_2h_0+h_4-c_2f(h_0)),\\
&& [x,y']= d_2x+d_3a'+h_3'\ (h_3'=h_3-d_1f(h_0)), \\
&& [x,h]= \lambda(h)y'+f_1(h)\ (f_1(h)=f(h)+\lambda(h)d_1h_0), \\
&& [y',h]=g(h),\quad [a',h]=h,\quad  [y',a']=h_1'\ (h_1'=h_1-c_2g(h_0)+d_1h_0),\\
&& \omega(y',a')=0,\quad \omega(x,a')=0,\quad \omega(x,y')=\omega(x,y)=1.
\end{eqnarray*}
For $[[x,a'],h]$, we have
 \begin{eqnarray*}
 [[x,a'],h] &=& c_1\lambda(h)y'+c_1f_1(h)+c_3h \\
 &=&[x,[a',h]]+[a',[h,x]]=\lambda(h)y'+\lambda(h)h_1'.
 \end{eqnarray*}
Then $c_1=1$ since $\lambda\not=0$, and $f_1(h)=-c_3h$ by taking $h\in \ker \lambda.$ For $[[x,y'],h]$, we have
\begin{eqnarray*}
[[x,y'],h]&=&[d_2x+d_3a'+h_3',h] =d_2\lambda(h)y'+d_2f_1(h)+d_3h \\
      &=&[x,[y',h]]+[y',[h,x]]+\omega (x,y')h=\lambda(g(h))y'+f_1(g(h))-g(f_1(h))+h.
\end{eqnarray*}
Take $h\in \ker \lambda$. Then $f_1(h)=-c_3h$ and $g(h)=0$. Considering $[[x,y'],h]|_H$, we have $$d_2c_3=d_3-1.$$
Finally, for $[[x,y'],a']$, we have
\begin{eqnarray*}
[[x,y'],a']&=&[d_2x+d_3a'+h_3',a']=d_2x+d_2c_3a'+d_2h_4'-h_3' \\
              &=&[x,[y',a']]+[y',[a',x]]+\omega (x,y')a' \\
              &=&\lambda(h_1')y'+d_2x+(d_3+1)a'+f_1(h_1')+h_3'-c_3h_1'-g(h_4').
\end{eqnarray*}
It follows that $$d_2c_3=d_3+1.$$
It is impossible.

\subsection{The proof of Theorem~\ref{classification}.} In summary, we complete Theorem~\ref{classification}.

\section*{Acknowledgements}
Z. Chen was partially supported by NNSF of China (11931009 and 12131012).

\end{document}